\documentclass[a4paper,14pt]{article}
\usepackage[english,russian]{babel}
\usepackage[left=25mm,right=15mm,top=20mm,bottom=20mm]{geometry}

\usepackage{verbatim}
\usepackage{amsmath}
\usepackage{amsthm}
\usepackage{amssymb}
\usepackage{delarray}
\usepackage{mathrsfs}
\usepackage{graphicx}
\usepackage{tikz}
\usepackage{fancyhdr}
\theoremstyle{plain}

\theoremstyle{definition}

\theoremstyle{remark}

 \allowdisplaybreaks
\usepackage{titlesec}
\titleformat{\section}
{\normalfont\bfseries}{}{1.2em}{}

\numberwithin{equation}{section}

\begin{document}

\begin{center}
\textbf{Recovering the characteristic functions of the Sturm-Liouville differential operators with singular potentials on star-type graph with cycle}

{\bf Vasilev S.~V.}
\end{center}

{\bf Abstract. } We consider Sturm-Liouville operators with singular potentials from the class $W^{-1}_2$ on star-type graph with cycle, which consist the edges with commensurable lengths. Asymptotic representation for eigenvalues for such operators is obtained. Recovering of the characteristic function the Sturm-Liouville operators with the singular potentials is considered. 
 


\section{1. Introduction.} \

This paper is devoted to the inverse spectral problems for differential operators on geometrical graphs. Differential operators on graphs are intensively studied by mathematicians in recent years and have applications in different branches of science and engineering. Inverse problem consists in recovering the potential on a graph from the given spectral characteristics. Recovering of the characteristic function of the Sturm-Liouville operators is the important part of solving such problem. The greatest success in the inverse spectral theory has been achieved for the classical Sturm-Liouville operator and afterwards for higher order differential operators [1-3].

The Sturm-Liouville operator with singular potential has been studied in papers [9-10]. The inverse problem on a finite interval were extensively studied in [11]. There are only few results for Sturm-Liouville operators with the singular potentials on graphs [27-29]. In this paper we consider the recovering of the characteristic function the Sturm-Liouville operators with the singular potentials on star-type graph with cycle, which consist the edges with commensurable lengths. Also we obtain the asymptotic representation for eigenvalues for such operators.

Let $G$ be a graph with of set of vertices $\{v_j\}_{j=0}^p$ and set of edges $\{e_j\}_{j=0}^p$, $e_j=[v_0, v_j]$, and edge $e_0$ generates the cycle. We suppose that the length of edge $e_j$ is equal to $|e_j|$. We consider each edge $e_j$ as a segment $[0, |e_j|]$ and parameterize it by the parameter $x\in [0, |e_j|]$. It is convenient for us to choose the orientation such that $x_j=|e_j|$ corresponds to the vertex $v_0$. Lengths of the edges $e_j$, $j=\overline{0,p}$, we consider as commensurable quantities.

A function $y$ on graph $G$ is considered as $y=[y_j(x_j)]_{j=0}^p$, where $y_j(x_j)$ are correspond to $e_j$, $x_j\in[0,|e_j|]$. Let $q=[q_j(x_j)]_{j=0}^p$ be a real-valued function on $G$ such that $q_j\in W_2^{-1}[0,|e_j|]$, i.e. $q_j(x_j)=\sigma'_j(x_j)$, where the derivative is considered in the sense of distributions. Function $\sigma=[\sigma_j(x_j)]_{j=0}^p$ we call the potential. The Sturm-Liouville differential operator on the edges $e_j$, $j=\overline{0,p}$ is defined by the following expression:
$$
\ell_jy_j:=-\big(y_j^{[1]}\big)'-\sigma_j(x_j)y_j^{[1]}-\sigma^2_j(x_j)y_j,
$$
where $y_j^{[1]}:=y_j'-\sigma_j(x)y_j$ - is a quasi-derivative, and
$$
dom(\ell_j)=\{y_j \ | \ y_j\in W^1_2[0,|e_j|], \ y_j^{[1]}\in W_1^1[0,|e_j|], \  \ell_jy_j\in L_2[0,|e_j|]\}.
$$
We consider the Sturm-Liouville equation on $G$:
\begin{equation}\label{1.1}
(\ell_jy_j)(x_j)=\lambda y_j(x_j), \ x_j\in(0,|e_j|),\ y_j \in dom(\ell_j),\ j=\overline{0,p}
\end{equation}
At the internal vertex $v_0$ we consider the following matching conditions:
\begin{equation}\label{1.2}
y_0(0)=y_k(|e_k|), \ k=\overline{0,p}, \quad \sum_{j=0}^py_j^{[1]}(|e_j|)=y_0^{[1]}(0).
\end{equation}
Let us consider the boundary value problem $L(G)$ for equation \eqref{1.1} with the matching conditions \eqref{1.2} and boundary conditions
\begin{equation}\label{1.3}
y_j^{[1]}(0)=0, \quad j=\overline{1,p}.
\end{equation}
Let us define by $\Lambda$ the eigenvalues of $L(G)$. We also consider the boundary value problem $L_k(G)$, $k=\overline{1,p}$, for equation \eqref{1.1} with the matching conditions \eqref{1.2} and boundary conditions
\begin{equation}\label{1.4}
y^{[1]}_j(0)=0,\quad j=\overline{1,p}\textbackslash k, \quad y_k(0)=0.
\end{equation}
The eigenvalues of $L_k(G)$ we define as $\Lambda_k$.

\section{2. Auxiliary propositions.} \

Let $C_j(x_j,\lambda),$ $S_j(x_j,\lambda)$  be the solutions of equation \eqref{1.1} on $e_j$, $j=\overline{0,p}$ under initial conditions
\begin{equation}\label{2.1}
C_j(0,\lambda)=S^{[1]}_j(0,\lambda)=1, \ \ C_j^{[1]}(0,\lambda)=S_j(0,\lambda)=0,
\end{equation}
From Liouville's formula we assume that $\left\langle C_j(x_j,\lambda),S_j(x_j,\lambda)\right\rangle=1$, where wronskian $\left\langle y,z\right\rangle=yz^{[1]}-y^{[1]}z$. The characteristic functions of the boundary value problem $L(G)$ we denote as $\Delta(\lambda,L(G))$. As in the classical case [16] one can show that the functions $M_j(\lambda)$, $j=\overline{1,p}$ are meromorphic in $\lambda$; namely:
\begin{equation}\label{2.3}
M_j(\lambda)=\frac{\Delta(\lambda, L_j(G))}{\Delta(\lambda,L(G))}.
\end{equation}
The vertex $v_0$ separate $G$ on two parts: $G=G_0\cup T$, where $G_0$ is the star-type graph with the internal vertex $v_0$ and $T$ is the cycle, generated by the edge $e_0$. Taking into account \eqref{1.2} è \eqref{1.3}, one can show, that
\begin{equation}\label{T}
\Delta(\lambda, L(T)) = C_0(|e_0|,\lambda)+S^{[1]}_0(|e_0|,\lambda)-2
\end{equation}
Denote by $B_r$, a the class of Paley-Wiener functions of exponential type not greater than $r \in \mathbb{R}$, belonging to $L_2(\mathbb R)$. It follows from [9-11] that
\begin{equation}\label{2.28}
C_j(|e_j|, \lambda) = \cos\rho|e_j| + \zeta_{j,e}(\rho), \quad S_j(|e_j|, \lambda) = \frac{\sin\rho|e_j|}{\rho} + \frac1\rho\zeta_{j,o}(\rho),
\end{equation}
where $\zeta_{j,e}(\rho)\in B_{|e_j|}$ are even functions and $\zeta_{j,o}(\rho)\in B_{|e_j|}$ are odd function. Clearly, that
$$
\zeta_{j,o}(\rho) = \int\limits_0^{|e_j|} K_{o}(t)\sin\rho tdt, \quad \zeta_{j,e}(\rho) = \int\limits_0^{|e_j|} K_{e}(t)\cos\rho tdt, \quad K_o, \ K_e \in L_2(0,|e_j|).
$$
Analogously [33], one can prove the following lemma.

{\bf Lemma 2.1. }{\it The characteristic functions of the boundary value problems $L(G)$ è $L_k(G)$, $k=\overline{1,p}$, admits the following representation
\begin{equation}\label{2.16}
\begin{array}{c}
\displaystyle \Delta(\lambda, L(G))=S_0(|e_0|,\lambda)\Delta(\lambda, L(G_0))+(-1)^p \Delta(\lambda, L(T))\prod\limits_{k=1}^pC_k(|e_k|, \lambda)\\[3mm]
\displaystyle \Delta(\lambda, L_j(G))=S_0(|e_0|,\lambda)\Delta_j(\lambda, L(G_0))+(-1)^p\Delta(\lambda, L(T))S_j(|e_j|,\lambda)\prod\limits_{k=\overline{1,p}\textbackslash\{j\}}C_k(|e_k|, \lambda)
\end{array}
\end{equation}
where
$$
\Delta(\lambda, L(G_0))=(-1)^p\sum\limits_{j=1}^pC_j^{[1]}(|e_j|,\lambda)\prod\limits_{i=1, \ i\ne j}^pC_i(|e_i|,\lambda)
$$
$$\begin{array}{c}
\displaystyle\Delta(\lambda, L_k(G_0))=(-1)^{p+1}\Big[S_k(|e_k|,\lambda))\sum\limits_{j=1, \ i\ne k} ^p C_j^{[1]}(|e_j|,\lambda)\prod\limits_{i=1, \ i\ne k,j}^pC_i(|e_i|,\lambda) + \\[3mm] \displaystyle S^{[1]}_k(|e_k|,\lambda))\prod\limits_{i=1}^pC_j(|e_j|,\lambda)\Big]
\end{array}
$$
}

\section{3. Asymptotic of the spectrums. Recovering the characteristic function}\

The eigenvalues $\Lambda$ can be numbered as $\{\lambda_{nk}\}_{k=\overline{0,\mu_0-1}, n\in\mathbb N} \cup \{\lambda_{nk}\}_{k=\overline{\mu_0, m}, n\in\mathbb Z}$, where $m \in \mathbb N$, $\mu_0$ is the multiplicity of the zero eigenvalue of the boundary value problem $L(G)$ with the zero potential.

{\bf Lemma 3.1. }{\it
Following asymptotic behavior are valid for the eigenvalues of the boundary value problem $L(G)$:
\begin{equation}\label{2.29}
\sqrt{\lambda_{nk}} =: \rho_{nk} = \left\{
\begin{array}{l}
\displaystyle \tau n + \varepsilon_{nk}, \quad k=\overline{0,\mu_0-1}, \quad n\in\mathbb N,\\[3mm]
\displaystyle |\tau n + \alpha_k| + \varepsilon_{nk}, \quad k=\overline{\mu_0, m}, \quad n\in\mathbb Z,
\end{array}
\right.
\end{equation}
where $\varepsilon_{nk}\in l_{2\mu_k}, \ \tau \in \mathbb R, \ \mu_k\in\mathbb N$.
}

{\bf Proof.} From \eqref{2.16} and \eqref{2.28} we assume
\begin{equation}\label{2.30}
\Delta(\lambda, L(G))=\Delta_0(\lambda, L(G)) + \zeta_{G}(\rho),
\end{equation}
where $\zeta_{G}(\rho) \in B_{|G|}$ are even and
$$
\Delta_0(\lambda, L(G)) = (-1)^p\Big[2\Big(\cos\rho|e_0|-1\Big)\prod\limits_{k=1}^p\cos\rho|e_k|
-\sin\rho|e_0|\sum\limits_{j=1}^p\sin\rho|e_j|\prod\limits_{i=1, \ i\ne
j}^p\cos\rho|e_i|\Big]
$$
Function $\Delta_0(\lambda, L(G))$ are the characteristic function of the boundary value problem $L(G)$ with zero potential. Define $d_0(\rho):=\Delta_0(\rho^2, L(G))$ and $d(\rho):=\Delta(\rho^2, L(G))$. From the commensurability of the lengths of the edges it follows, that function $d_0(\rho)$ are periodic. Smallest period $\tau>0$ of the function $d_0(\rho)$ is in existence (see [23]). Its sufficient to investigate zeros of the function $d_0(\rho)$ on the interval $[0,\tau/2]$, because function $d_0(\rho)$ is even. Define zeros of the function $d_0(\rho)$ on $[0,\tau/2]$ as (counting with the multiplicities) $\alpha_k$, $k=\overline{0,m}$, where $m \in \mathbb N$. We define the multiplicity of the  $\alpha_k$ as $\mu_k$, $k=\overline{0,m}$. Clearly, that $\alpha_0=...=\alpha_{\mu_0}=0$. For simplicity we assume, that $\tau/2$ are not the zero of $d_0(\rho)$. Then we conclude:

\begin{equation}\label{2.31}
\rho^{0}_{nk} := \left\{
\begin{array}{l}
\displaystyle \tau n, \quad k=\overline{0,\mu_0-1}, \quad n\in\mathbb N,\\[3mm]
\displaystyle |\tau n + \alpha_k|, \quad k=\overline{\mu_0, m}, \quad n\in\mathbb Z
\end{array}
\right.
\end{equation}
Let us consider $\delta:=\min\limits_{\alpha_k\ne\alpha_j, \ k,j=\overline{0,m}}|\alpha_k-\alpha_j|$. By the sequence $\rho^{0}_n$, $n\in\mathbb N$ we assume the sequence $\rho^{0}_{nk}$, $(n,k)\in \{\overline{0,\mu_0-1}\} \times \mathbb N \bigcup \{\overline{\mu_0, m}\} \times \mathbb Z$ put in order of increasing, such that $\rho^{0}_n\ne \rho^{0}_s$ for $n\ne s$. Define $\Gamma_n := \{\lambda, |\lambda| = (\rho^{0}_n+\delta)^2\}$, $n\in \mathbb N$. Clearly, that
$$
|\zeta_{|G|, e}(\rho)| \le C_1 e^{|G||\rho^{0}_n +\delta|}, \quad |\Delta_0(\lambda, L(G))|\ge C_2 e^{|G||\rho^{0}_n +\delta|}, \quad \rho \in \partial\Gamma_n
$$
Analogously [7], by the Rouche's theorem, one can show:
\begin{equation}\label{2.33}
\rho_{nk} := \left\{
\begin{array}{l}
\displaystyle \tau n + \varepsilon_{nk}, \quad k=\overline{0,\mu_0-1}, \quad n\in\mathbb N,\\[3mm]
\displaystyle |\tau n + \alpha_k| + \varepsilon_{nk}, \quad k=\overline{\mu_0, m}, \quad n\in\mathbb Z,
\end{array}
\right.
\quad \varepsilon_{nk} = o(1), \quad n \to\infty.
\end{equation}
Without restricting the generality we consider the case $n\ge 0$ with $k=\overline{\mu_0, m}$. Substituting \eqref{2.31} into \eqref{2.30}, we obtain
\begin{equation}\label{analityc}
d_0(\alpha_k + \varepsilon_{nk}) + \zeta_{|G|, e}(\rho_{nk}) = 0
\end{equation}
Clearly, that
\begin{equation}\label{2.kdz}
\zeta_{G, e}(\rho_{nk}) = \int\limits_0^{|G|} K(t)\cos(n\tau + \alpha_k + \varepsilon_{nk}) tdt, \quad K(t) \in L_2(0,|G|)
\end{equation}
Let $N^*$ be a such number, that $\forall n > N^*$ $\varepsilon_{nk} < 1/4$. Then from the kadec-1/4 theorem in the complex case (see [24]) we conclude, that $\{\cos(n+\varepsilon^*_{nk})x\}_n$ and $\{\sin(n+\varepsilon^*_{nk})x\}_n$, where
$$
\varepsilon^*_{nk}=
\left\{
\begin{array}{c}
0, \quad n \le N^*, \\[3mm]
\varepsilon_{nk}/\tau, \quad n > N^*,
\end{array}
\right.
$$
are the Riesz basis. Then change the variables $s=\tau t$ and taking into account
$$
\begin{array}{c}
\displaystyle\cos(n\tau + \alpha_k + \varepsilon_{nk})t = \cos\alpha_k t\cos(n\tau + \varepsilon_{nk})t-\sin\alpha_k t\sin(n\tau + \varepsilon_{nk})t,
\end{array}
$$
from \eqref{2.kdz} we assume, that
\begin{equation}\label{2.kdz1}
\zeta_{G, e}(\rho_{nk}) = \int\limits_0^{\pi} A_k(t)\cos(n + \varepsilon^*_{nk}) tdt + \int\limits_0^{\pi} B_k(t)\sin(n + \varepsilon^*_{nk}) tdt \in l_2, \quad A_k, B_k \in L_2(0,\pi)
\end{equation}
By virtue of the fact, that function $d_0(\rho)$ are analytic at a point $\alpha_k$, we obtain, that
$$
d_0(\alpha_k+\varepsilon_{nk}) = \sum\limits_{j=0}^{\mu_k}d^{(j)}_0(\alpha_k)\frac{\varepsilon_{nk}^j}{j!} + O(\varepsilon_{nk}^{\mu_k+1}), \quad n\to\infty
$$
From \eqref{analityc} it follows, that
$$
\frac{1}{\mu_k!}d_0^{(\mu_k)}(\alpha_k)\varepsilon_{nk}^{\mu_k} + O(\varepsilon_{nk}^{\mu_k+1}) + \zeta_{|G|, e}(\rho_{nk}) = 0, \quad n \to \infty
$$
Consequently, 
$$
\varepsilon^{\mu_k}_{nk} = \zeta_{|G|, e}(\rho_{nk})\frac{1}{1+O(\varepsilon_{nk})} \in l_2, \quad \varepsilon_{nk} \in l_{2\mu_k}.
$$
Thus, we obtain $\eqref{2.29}$.$\hfill\Box$

Define
$$
\lambda^1_{nk} :=
\left\{
\begin{array}{c}
\lambda_{nk}, \quad \lambda_{nk} \ne 0, \\[3mm]
1, \quad \lambda_{nk} = 0,
\end{array}
\right.
\quad
\lambda^{01}_{nk} :=
\left\{
\begin{array}{c}
\lambda^0_{nk}, \quad \lambda^0_{nk} \ne 0, \\[3mm]
1, \quad \lambda^0_{nk} = 0,
\end{array}
\right.
\quad
\lambda^{01}_{jnk} :=
\left\{
\begin{array}{c}
\lambda^0_{jnk}, \quad \lambda^0_{jnk} \ne 0, \\[3mm]
1, \quad \lambda^0_{jnk} = 0,
\end{array}
\right.
$$
The following theorem give us the formulas for the recovering the characteristic function of the boundary value problem $L(G)$ and $L_j(G)$, $j=\overline{1,p}$, by the spectra $\Lambda$ è $\Lambda_j$, $j=\overline{1,p}$ respectively.

{\bf Theorem 3.1. }{\it The specification of the spectrums $\Lambda$ and $\Lambda_j$ uniquely determines the characteristic functions respectively by the formula
\begin{equation}\label{2.39}
\begin{array}{c}
\displaystyle\Delta(\lambda, L(G)) =
(-1)^{\mu_0}\frac{\partial^{\mu_0}}{\partial\lambda^{\mu_0}}\Delta_0(\lambda, L(G)) \Big|_{\lambda = 0} \prod\limits_{n=0}^{\infty}\prod\limits_{k=0}^{\mu_0-1}\frac{\lambda_{nk}-\lambda}{\lambda^{01}_{nk}}\prod\limits_{n=-\infty}^{\infty}\prod\limits_{k=\mu_0}^{m}\frac{\lambda_{nk}-\lambda}{\lambda^{01}_{nk}}.\\[3mm]
\displaystyle \Delta(\lambda, L_j(G)) =
(-1)^{\mu_0}\frac{\partial^{\mu_0}}{\partial\lambda^{\mu_0}}\Delta_0(\lambda, L_j(G)) \Big|_{\lambda = 0} \prod\limits_{n=0}^{\infty}\prod\limits_{k=0}^{\mu_0-1}\frac{\lambda_{jnk}-\lambda}{\lambda^{01}_{jnk}}\prod\limits_{n=-\infty}^{\infty}\prod\limits_{k=\mu_0}^{m}\frac{\lambda_{jnk}-\lambda}{\lambda^{01}_{jnk}}.
\end{array}
\end{equation}
}

{\bf Proof.} Let us consider the case of the boundary value problem $L(G)$. Using Hadamard's factorization theorem, one can show
\begin{equation}\label{2.40}
\frac{\Delta(\lambda, L(G))}{\Delta_0(\lambda, L(G))} = \frac{C}{C_0}
\prod\limits_{n=0}^{\infty}\prod\limits_{k=0}^{\mu_0-1}\frac{\lambda^{01}_{nk}}{\lambda^1_{nk}}\Big(1 +\frac{\lambda_{nk} - \lambda_{nk}^0}{\lambda^0_{nk} - \lambda}\Big) \prod\limits_{n=-\infty}^{\infty}\prod\limits_{k=\mu_0}^{m}\frac{\lambda^{01}_{nk}}{\lambda^1_{nk}}\Big(1 +\frac{\lambda_{nk} - \lambda_{nk}^0}{\lambda^0_{nk} - \lambda}\Big),
\end{equation}
Clearly, that $\forall \lambda \in [-\infty, 0]$
$$
\Big|\sum\limits_{k=0}^{\mu_0-1}\frac{\lambda_{nk} - \lambda_{nk}^0}{\lambda^0_{nk} - \lambda}\Big| \le \sum\limits_{k=0}^{\mu_0-1}\Big|\frac{\lambda_{nk} - \lambda_{nk}^0}{\lambda^0_{nk}}\Big| \le \sum\limits_{k=0}^{\mu_0-1}\Big|\frac{2\varepsilon_{nk}}{\rho^0_{nk}}\Big| + \sum\limits_{k=0}^{\mu_0-1}\Big|\frac{\varepsilon_{nk}}{\rho^0_{nk}}\Big|^2,
$$
Using H{\"o}lder's inequality and $\varepsilon_{nk} \in l_{2\mu_k}$, one can show
$$
\sum\limits_{n=0}^{\infty} \sum\limits_{k=0}^{\mu_0-1}\Big|\frac{2\varepsilon_{nk}}{\rho^0_{nk}}\Big| + \sum\limits_{n=0}^{\infty} \sum\limits_{k=0}^{\mu_0-1}\Big|\frac{\varepsilon_{nk}}{\rho^0_{nk}}\Big|^2 < \infty,
$$
Then using Weierstrass M-test, we obtain the uniform convergence of the series
$$
\sum\limits_{n=0}^{\infty}\sum\limits_{k=0}^{\mu_0-1}\frac{\lambda_{nk} - \lambda_{nk}^0}{\lambda^0_{nk} - \lambda},
$$
where $\lambda \in [-\infty, 0]$. Then we obtain the uniform convergence of the infinite product
$$
\prod\limits_{n=0}^{\infty}\prod\limits_{k=0}^{\mu_0-1}\Big(1 +\frac{\lambda_{nk} - \lambda_{nk}^0}{\lambda^0_{nk} - \lambda}\Big)
$$
Clearly, that
$$
\lim\limits_{\lambda\to -\infty}\prod\limits_{n=0}^{N}\prod\limits_{k=0}^{\mu_0-1}\Big(1 +\frac{\lambda_{nk} - \lambda_{nk}^0}{\lambda^0_{nk} - \lambda}\Big) = 1, \quad \forall N \in \mathbb N.
$$
It follows from the Moore-Osgood theorem, that
\begin{equation}\label{2.42}
\lim\limits_{\lambda\to -\infty} \lim\limits_{N \to\infty}\prod\limits_{n=0}^{N}\prod\limits_{k=0}^{\mu_0-1}\Big(1 +\frac{\lambda_{nk} - \lambda_{nk}^0}{\lambda^0_{nk} - \lambda}\Big) =
\lim\limits_{N \to\infty}\lim\limits_{\lambda\to -\infty} \prod\limits_{n=0}^{N}\prod\limits_{k=0}^{\mu_0-1}\Big(1 +\frac{\lambda_{nk} - \lambda_{nk}^0}{\lambda^0_{nk} - \lambda}\Big) =1
\end{equation}
Using $\varepsilon_{nk}\in l_{2\mu_k}$, one can show the convergence of the infinite product
$$
\prod\limits_{n=0}^{\infty}\prod\limits_{k=0}^{\mu_0-1}\frac{\lambda^{1}_{nk}}{\lambda^{01}_{nk}} = \prod\limits_{n=0}^{\infty}\prod\limits_{k=0}^{\mu_0-1}\Big(1+2\frac{\varepsilon_{nk}}{\rho_{nk}^0}+ \Big(\frac{\varepsilon_{nk}}{\rho_{nk}^0}\Big)^2\Big)
$$
Clearly, that
$$
C_0 = (-1)^{\mu_0}\frac{\partial^{\mu_0}}{\partial\lambda^{\mu_0}}\Delta_0(\lambda, L(G)) \Big|_{\lambda = 0},
$$
Using lemma 2.1 and \eqref{2.30}, we assume, that $\lim\limits_{\lambda \to -\infty}\Delta(\lambda, L(G))/\Delta_0(\lambda, L(G)) = 1$. Consequently, using \eqref{2.40} we obtain \eqref{2.39}.
Analogously, one can prove \eqref{2.39} for the boundary value problem $L_j(G)$. $\hfill\Box$

{\bf Acknowledgment. } This work was supported in part by Grants \textnumero 16-01-00015,  \textnumero 17-51-53180 of the Russian Foundation for Basic Research, by Grant \textnumero 1.1660.2017/4.6 of the Russian Ministry of Education and Science.

\section{References} \

1. V.A. Marchenko, Sturm-Liouville Operators and Their Applications, Naukova Dumka, Kiev, 1977; English transl., Birkhauser, 1986.

2. B.M. Levitan, Inverse Sturm-Liouville Problems, Nauka, Moscow, 1984; English transl., VNU Sci. Press, Utrecht, 1987.

3. G. Freiling and V.A. Yurko, Inverse Sturm-Liouville Problems and their Applications,NOVA Science Publishers, New York, 2001.

4. Å. Trubowitz, Inverse Spectral Theory, New York, Academic Press, 1987.

5. K.Chadan, D. Colton, L. Paivarinta and W. Rundell, An introduction to inverse scattering and inverse spectral problems, SIAM Monographs on Mathematical Modeling and Computation.SIAM, Philadelphia, PA, 1997.

6. R. Beals , P. Deift and C. Tomei, Direct and Inverse Scattering on the Line, Math. Surveys and Monographs, v.28. Amer. Math. Soc. Providence: RI, 1988.

7. G. Freiling, V.A. Yurko, Inverse Sturm-Liouville Problems and Their Applications, NOVA Science Publishers, New York, 2001.

8. V.A. Yurko, Method of Spectral Mappings in the Inverse Problem Theory, Inverse and Illposed Problems Series. VSP, Utrecht, 2002.

9. R.O. Hryniv and Ya.V. Mykytyuk, Inverse spectral problems for Sturm-Liouville operators with singular potentials, Inverse Problems 19, 2003.

10. R.O. Hryniv and Ya.V. Mykytyuk, Transformation operators for Sturm-Liouville operators with singular potentials, Mathematical Physics, Analysis and Geometry 7, 2004.

11. A. M. Savchuk, A. A. Shkalikov, Sturm-Liouville operators with singular potentials, Math. Notes 66:6, 1999.

13. M.I. Belishev, Boundary spectral inverse problem on a class of graphs (trees) by the BC method, Inverse Problems 20 (2004), 647-672pp.

14. V.A. Yurko, Inverse spectral problems for Sturm-Liouville operators on graphs, Inverse Problems 21, 2005.

15. G. Freiling and V.A. Yurko, Inverse problems for differential operators on trees with general matching conditions, Applicable Analysis 86, no.6 (2007), 653-667pp.

16. B.M. Levitan and I.S. Sargsyan, Introduction to Spectral Theory, AMS Transl. of Math. Monogr. 39, Providence, 1975, 525p.

17. R. Bellmann and Coike, Differential-difference Equations, 1963, 548p.

18. Freiling G., Yurko V.A. Inverse nodal problems for differential operators on graphs with a cycle. Tamkang J. Math. 41, N1. 2010, 50-54pp.

19. Yurko V.A. Inverse problems for Sturm-Liouville operators on bush-type graphs. Inverse Problems 25, N10. 2009, 125-127pp.

20. Yurko V.A. Uniqueness of recovering Sturm-Liouville operators on A-graphs from spectra. Results in Mathematics 55, N1-2. 2009, 199-207pp.

21. M.Y. Ignatiev. Inverse spectral problem for Sturm-Liouville operator on non-compact A-graph. Uniqueness result. Department of Mathematics, Saratov University, 25p.

22. M.Y. Ignatiev, Inverse scattering problem for Sturm-Liouville operator on one-vertex noncompact graph with a cycle, Tamkan J. of Mathematics 42, N3.011, 154-166pp.

23. J.B. Conway, Functions of One Complex Variable, vol.I, 2nd edn., Springer--Verlag, New York, 1995, 412p.

23. A. Yu. Evnin, Polynomial as a sum of periodic functions, Vestn. Yuzhno-Ural. Gos. Un-ta. Ser. Matem. Mekh. Fiz., 2013.

24. P. Vellucci, A simple pointview for kadec-1/4 theorem in the complex case, Ricerche di Matematica, 2014.

25. Taylor, Angus E. (2012). General Theory of Functions and Integration. Dover Books on Mathematics Series. p. 140.

26. Yurko, V. A. Inverse problems for Sturm-Liouville operators on graphs with a cycle, Operators and Matrices 2:4, 2008.

27. Bondarenko N.P. A 2-edge partial inverse problem for the Sturm-Liouville operators with singular potentials on a star-shaped graph, Tamkang Journal of Mathematics (2018), Vol. 49, No. 1. Pp. 49-66.

28. Yang C.-F., Bondarenko N.P. A partial inverse problem for the Sturm-Liouville operator on the lasso-graph, Inverse Problems and Imaging (2019), Vol. 13, No. 1.

29. G. Freiling, M. Ignatiev, V. Yurko, An inverse spectral problem for Sturm-Liouville operators with singular potentials on star-type graphs, 2008, 397-409pp.

\end{document}